# USE OF AUXILIARY INFORMATION IN VARIANCE ESTIMATION

[1]Jayant Singh, [2]Viplav K. Singh, [2]Sachin Malik and †[2]Rajesh Singh

[1]Department of Statistics, Rajasthan University, Jaipur

[2]Department of Statistics, BHU

†Corresponding author (rsinghstat@gmail.com)

**Abstract**

This paper proposes a class of ratio type estimators of finite population variance, when the population variance of an auxiliary character is known. Asymptotic expression for mean square error (MSE) is derived and compared with the mean square errors of some existing estimators. An empirical study is carried out to illustrate the performance of the constructed estimator over others.

## 1. INTRODUCTION

In manufacturing industries and pharmaceutical laboratories sometimes researchers are interested in the variation of their products. To measure the variations within the values of study variable y, the problem of estimating the population variance of $S_y^2$ of study variable y received a considerable attention of the statistician in survey sampling including Isaki (1983), Singh and Singh (2001, 2003), Jhajj et al. (2005), Kadliar and Cingi (2007), Singh et al. (2008), Grover (2010), Singh et al. (2011) and Singh and Solanki (2012) have suggested improved estimator for estimation of $S_y^2$.

Let us consider a finite population $U = (U_1, U_2, U_2........U_N)$ having N units and let y and x are the study and auxiliary variable with means $\overline{X}$ and $\overline{Y}$ respectively. Let us suppose that a sample of size n is drawn from the population using simple random sampling without replacement (SRSWOR) method. Let $s_y^2 = \sum_{i=1}^{n}(y_i - \overline{y})^2 /(n-1)$ and

$s_x^2 = \sum_{i=1}^{n}(x_i - \bar{x})^2/(n-1)$ be the sample variances for variables y and x, which are unbiased estimators for $S_Y^2 (= \sum_{i=1}^{N}(y_i - \bar{Y})^2/(N-1))$ and $S_x^2 (= \sum_{i=1}^{N}(x_i - \bar{X})^2/(N-1))$, respectively. Let

$$C_y = \frac{S_y}{\bar{Y}} \text{ and } C_x = \frac{S_x}{\bar{X}}.$$

Now, assume that the problem is to estimate the population variance $S_Y^2$ of the study variable y which is highly correlated with the auxiliary variable x.

In our paper, we use $\theta = \left(\frac{1}{n} - \frac{1}{N}\right)$.

Let,

$$\mu_{pq} = \frac{1}{N}\sum_{i=1}^{N}(y_i - \bar{Y})^p(x_i - \bar{X})^q,$$

$$\lambda_{22} = \frac{\mu_{pq}}{\mu_{20}^{p/2}\mu_{02}^{q/2}},$$

$$\beta_{2y} = \frac{\lambda_{40}}{\mu_{20}^2}, \quad \beta_{2y}^* = (\beta_{2y} - 1),$$

$$\beta_{2x} = \frac{\lambda_{40}}{\mu_{20}^2}, \quad \beta_{2x}^* = (\beta_{2x} - 1),$$

$$\lambda_{22}^* = \lambda_{22} - 1.$$

Let $e_0 = \left\{\frac{s_y^2}{S_y^2} - 1\right\}$ and $e_1 = \left\{\frac{s_x^2}{S_x^2} - 1\right\}$

such that,

$E(e_0) = E(e_0) = 0,$

and up to the first order of approximation,

$E(e_0^2) = \theta\beta_{2y}^*,$
$E(e_1^2) = \theta\beta_{2x}^*,$
$E(e_0 e_1) = \theta\lambda_{22}^*.$

In this paper, under SRSWOR, we have suggested a general family of estimators for estimating the population variance $S_y^2$. The expression of MSE, up to the first order of approximation have been obtained.

2. **Existing estimators**

The variance of the usual unbiased estimator $\hat{S}_y^2$ is given by

$$\text{var}(\hat{S}_y^2) = \theta\, S_y^4 \beta_{2y}^* \tag{2.1}$$

Isaki (1983) suggested following ratio and regression estimators

$$\hat{S}_R^2 = s_y^2 \left( \frac{S_x^2}{s_x^2} \right) \tag{2.2}$$

and $\quad \hat{S}_{\text{Reg}}^2 = s_y^2 + b^*(S_x^2 - s_x^2) \tag{2.3}$

The MSE's of the ratio and regression estimators are respectively, given by

$$\text{MSE}(S_R^2) = \theta S_y^4 \left( \beta_{2y}^* + \beta_{2x}^* - 2\lambda_{22}^* \right) \tag{2.4}$$

$$\text{MSE}(t_3) = \theta S_y^4 \beta_{2y}^* \left( 1 - \frac{\lambda_{22}^{*2}}{\beta_{2y}^* \beta_{2x}^*} \right) \tag{2.5}$$

where,

$$\rho^*_{(S_y^2, S_x^2)} = \frac{\lambda_{22}^*}{\sqrt{\beta_{2y}^*}\sqrt{\beta_{2x}^*}}.$$

Using known values of some population parameter(s), Kadilar and Cingi (2006) suggested the following variance estimators –

$$\hat{S}^2_{KC1} = s_y^2 \left( \frac{S_X^2 + C_X}{s_X^2 + C_X} \right) \tag{2.6}$$

$$\hat{S}^2_{KC2} = s_y^2 \left( \frac{S_X^2 + \beta_{2x}}{s_X^2 + \beta_{2x}} \right) \tag{2.7}$$

$$\hat{S}^2_{KC3} = s_y^2 \left( \frac{S_X^2 \beta_{2x} + C_X}{s_X^2 \beta_{2x} + C_X} \right) \tag{2.8}$$

$$\hat{S}^2_{KC4} = s_y^2 \left( \frac{S_X^2 C_X + \beta_{2x}}{s_X^2 C_X + \beta_{2x}} \right) \tag{2.9}$$

The MSE of $\hat{S}_{KC_i}$ (i = 1,2,3,4) to the first order of approximation is given by-

$$\text{MSE}(\hat{S}^2_{KCi}) = \theta S_y^4 \left[ \beta^*_{2y} + p_i^2 \beta^*_{2x} - 2 p_i \lambda^*_{22} \right] \tag{2.10}$$

where,

$$p_1 = \frac{S_X^2}{S_X^2 + C_X}, \; p_2 = \frac{S_X^2}{S_X^2 + \beta_{2x}}, \; p_3 = \frac{S_X^2 \beta_{2x}}{S_X^2 \beta_{2x} + C_X}, \; p_4 = \frac{S_X^2 C_X}{S_X^2 C_X + \beta_{2x}}.$$

Kadilar and Cingi (2006) also suggested the following estimator

$$\hat{S}^2_{KC} = \alpha_1 s_y^2 + \alpha_2 \left( s_y^2 \frac{S_X^2}{s_X^2} \right) \tau \tag{2.11}$$

where $\alpha_1 + \alpha_2 = 1$ and $\tau = \frac{1 + \lambda c_{yx}}{1 + \lambda c_x^2}$.

MSE of $\hat{S}^2_{KC}$ to the first order of approximation is given by

$$\text{MSE}(\hat{S}^2_{KC}) = \theta S_y^4 \left[ z^2 \beta^*_{2y} + \alpha_2^{*2} \tau^2 \beta^*_{2x} - 2 \tau z \alpha_2^* \lambda^*_{22} \right] \tag{2.12}$$

where,

$$z = \alpha_1^* + \alpha_2^* \tau$$

$$\alpha_1^* = \frac{\beta^*_{2y}(\tau - 1) + \beta^*_{2x} \tau + (1 - 2\tau)\lambda^*_{22}}{\beta^*_{2y}\left\{ \frac{(1-\tau)^2}{\tau} \right\} + 2\lambda^*_{22}(1-\tau) + \beta^*_{2x}\tau}, \alpha_2^* = 1 - \alpha_1^*$$

Using Ray and Sahai (1980) estimator, Gupta and Shabbir (2007) proposed a hybrid class of estimators of $S_y^2$ as

$$S^{2(\alpha)}_{PR} = \left[ d_1 s_y^2 + d_2 (S_X^2 - s_X^2) \right] \left[ 2 - \left( \frac{s_X^2}{S_X^2} \right)^\alpha \right] \tag{2.13}$$

where α is a constant.

The MSE of the estimator $S_{PR}^{2(\alpha)}$ is given by

$$\text{MSE}(S_{PR}^{2(\alpha)}) = S_y^4 + d_1^2 S_y^4 A_1 + d_2^2 S_x^4 A_2 + 2d_1 d_2 S_y^2 S_x^2 A_3 - 2d_1 S_y^4 A_4 - 2d_2 S_y^2 S_x^2 A_5 \quad (2.14)$$

where,

$$A_1 = 1 + \theta(\beta_{2y}^* + \alpha \beta_{2x}^* - 4\alpha \lambda_{22}^*),$$

$$A_2 = \theta \beta_{2x}^*,$$

$$A_3 = \theta(2\alpha \beta_{2x}^* - \lambda_{22}^*),$$

$$A_4 = 1 - \alpha\theta \left( \lambda_{22}^* + \frac{(\alpha - 1)}{2} \beta_{2x}^* \right),$$

$$A_5 = \alpha\theta \beta_{2x}^*.$$

## 3. Proposed estimator

Motivated by Sahai and Ray (1980) and Singh and Solanki (2012), we have suggested the following generalized class of estimators for population variance $S_y^2$ of study variable y as

$$T_1 = w_1 s_y^2 \left\{ \frac{cS_X^2 - ds_x^2}{(c-d)S_X^2} \right\}^m + w_2 s_y^2 \left\{ 2 - \left( \frac{s_x^2}{S_X^2} \right)^w \right\} \quad (3.1)$$

Where ($w_1$, $w_2$) are constants to be determined such that MSE of generalized class of estimator T is minimum and m and w are constants used for generating different members of the class and c and d are either constants or function of known parameters of auxiliary variable x ( for choice of c and d refer to Singh and Kumar (2011)).

Expressing equation (3.1) in the form of $e_i$'s, (i=01) we have,

$$T = w_1(e_0 + 1)S_Y^2 \left[ \frac{cS_x^2 - d(e_1 + 1)S_x^2}{(c-d)S_x^2} \right]^m + w_2(e_1 + 1)S_Y^2 \left[ 1 - we_1 - \frac{w(w-1)}{2} e_1^2 \right] \quad (3.2)$$

Expanding the right hand side of expression (3.2), to the first order of approximation and subtracting $S_Y^2$ from both sides, we have

$$T - S_y^2 = S_y^2 \left\{ w_1 \left[ 1 - mAe_1 - m(m-1)\frac{A^2 e_1^2}{2} + e_0 - mAe_0 e_1 \right] \right.$$

$$\left. + w_2 \left[ 1 - we_1 - \frac{w(w-1)}{2} e_1^2 + e_0 - we_0 e_1 \right] - 1 \right\} \quad (3.3)$$

where, $A = \dfrac{d}{(c-d)}$.

Squaring both sides of equation (3.3) and then taking expectations, we get the MSE up to the first order of approximation of the estimator T, as

$$\text{MSE}(T) = S_y^4 \left\{ 1 + w_1^2 B_1 + w_2^2 B_2 + 2w_1 w_2 B_3 - 2w_1 B_4 - 2w_2 B_5 \right\} \quad (3.4)$$

where,

$$B_1 = 1 + \theta \left[ \beta_{2y}^* + mA^2 \beta_{2x}^* - 4mA\lambda_{22}^* \right]$$
$$B_2 = 1 + \theta \left[ \beta_{2y}^* + w\beta_{2x}^* - 4w\lambda_{22}^* \right]$$
$$B_3 = 1 + \theta \left[ \beta_{2x}^* \left( mwA - \frac{w(w-1)}{2} - \frac{m(m-1)}{2} \right) + \beta_{2y}^* - 2w\theta^* - 2mA\lambda_{22}^* \right]$$
$$B_4 = 1 - \theta \left[ \beta_{2x}^* A^2 \frac{m(m-1)}{2} + mA\lambda_{22}^* \right]$$
$$B_5 = 1 - \theta \left[ \beta_{2x}^* \frac{w(w-1)}{2} + w\lambda_{22}^* \right]$$

Minimizing expression (3.4) with respect to $w_1$ and $w_2$, we get the optimum values of $w_1$ and $w_2$ as

$$w_1(\text{opt}) = \frac{B_2 B_4 - B_3 B_5}{B_1 B_2 - B_3^2}$$
$$w_2(\text{opt}) = \frac{B_1 B_5 - B_3 B_4}{B_1 B_2 - B_3^2} \quad (3.5)$$

Using these optimum values of $w_1$ and $w_2$ in expression (3.4), we get the minimum MSE of the estimator T.

4. **Numerical illustration**

We use data in and Kadilar and Cingi (2007) to compare efficiencies between the traditional and proposed estimators in the simple random sampling.

So consider 104 villages of the East Anatolia Region in Turkey. Take the following variables:

y : level of apple production (1 unit = 100 tones)

x : number of apple trees (1 unit = 100 trees)

The values of required parameters of the population are:

$N = 104$, $S_y = 11.669964$, $S_x = 23029.072$, $C_y = 1.866$, $C_x = 1.653$, $\rho_{yx} = = 0.865$,

$C_{yx} == 2.668$, $\beta_{2(y)} = 16.523$, $\beta_{2(x)} = 17.516$, $\lambda_{22} = 14.398$.

Table 4.1 : MSE of the estimators

| Estimators | MSE |
|---|---|
| $S_y^2$ | 11627.2 |
| $S_R^2$ | 3927.166 |
| $S_{KC_1}^2$ | 3927.178 |
| $S_{KC_2}^2$ | 3927.178 |
| $S_{KC_3}^2$ | 3927.178 |
| $S_{KC_4}^2$ | 3927.178 |
| $S_{`KC}^2$ | 3473.024 |
| $S_{`Reg}^2$ | 3927.178 |
| $S_{`PR}^{2(\alpha=0)}$ | 2934.649 |
| $S_{PR}^{2(\alpha=1)}$ | 8721.148 |
| $S_{PR}^{2(\alpha=-1)}$ | 14832.09 |
| Proposed T(n=-1) | 347.6189 |
| T(n=0) | 7792.016 |
| T(n=1) | 11257.42 |

## 5. Conclusion

From theoretical discussion in section 3 and results of the numerical example, we infer that the proposed estimator 'T' under optimum condition performs better than usual estimator $S_y^2$, Isaki's (1983) estimator, Kadilar and Cingi's (2006) estimators and Gupta and Shabbir (2007) estimators.

## REFERENCES

Grover, L. K . (2010): A correction note on improvement in variance estimation using auxiliary information. Commun. Stat. Theo. Meth. 39:753–764.

Isaki, C.T.(1983): Variance estimation using auxiliary information. Journal of American Statistical Association 78, 117–123.

Jhajj, H .S., Sharma, M. K. and Grover, L. K. (2005) : An efficient class of chain estimators of population variance under sub-sampling scheme. J. Japan Stat. Soc., 35(2), 273-286.

Kadilar, C. and Cingi, H. (2006) : Improvement in variance estimation using auxiliary information. Hacettepe Journal of Mathematics and Statistics 35 (1), 111–115.

Kaur,(1985): An efficient regression type estimator in survey sampling, Biom. Journal 27 (1),107–110.

Ray, S.K. and Sahai, A (1980) : Efficient families of ratio and product type estimators, Biometrika 67(1), 211–215.

Shabbir, J. and Gupta, S. (2007): On improvement in variance estimation using auxiliary information. Commun. Statist. Theor. Meth. 36(12):2177–2185.

Singh, H.P. and Singh, R. (2001): Improved ratio-type estimator for variance using auxiliary information. J Indian Soc Agric Stat 54(3):276–287.

Singh, H.P. and Singh, R. (2003): Estimation of variance through regression approach in two phase sampling. Aligarh Journal of Statistics, 23, 13-30.

Singh, R., Chauhan, P., Sawan, N. and Smarandache, F.(2008): Almost unbiased ratio and product type estimator of finite population variance using the knowledge of kurtosis of an auxiliary variable in sample surveys. Octogon Mathematical Journal, Vol. 16, No. 1, 123-130.

Singh, R. , Kumar, M., Singh, A.K. and Smarandache, F. (2011): A family of estimators of population variance using information on auxiliary attribute. Studies in sampling techniques and time series analysis. Zip publishing, USA. (Singh, R. and Smarandache, F. editor).

Singh, H. P. and Solanki, R.S. (2012): A new procedure for variance estimation in simple random sampling using auxiliary information. Stat. Paps. DOI 10.1007/s00362-012-0445-2.